\documentclass[a4,12pt]{article}
\usepackage{amsmath}
\usepackage{latexsym}
\author{Magnus Fontes}
\newtheorem{definition}{Definition}
\newtheorem{corollary}{Corollary}
\newtheorem{theorem}{Theorem}
\newtheorem{lemma}{Lemma}
\numberwithin{theorem}{section}
\numberwithin{equation}{section}
\numberwithin{definition}{section}
\numberwithin{corollary}{section}
\numberwithin{lemma}{section}
\title{Initial-Boundary Value Problems for Parabolic Equations.}
\date{}
\begin{document}
\maketitle

\section{Introduction.}
In this paper we prove  new  existence and uniqueness results for 
weak solutions to 
non-homogeneous initial-boundary value problems 
for parabolic equations of the form
\begin{subequations}\label{introeq}
\begin{align}
 \frac{\partial u}{\partial t} - \nabla_x \cdot A(x,t,\nabla_x u)
 &=f \quad
\mbox{in $\mathcal{D}'(Q_+)$}\\
u&=g \quad \mbox{on $ (\Omega \times \{0\}) \cup (\partial \Omega \times {\bf R}_+)$}.
\end{align}
\end{subequations}
Here $\Omega$ is an open and bounded set 
in ${\bf R}^n$ and  $Q_+ = \Omega \times {\bf R_+}$.
Precise structural conditions for $A(\cdot,\cdot,\cdot)$ are given in
Section 4, but the 
model is the following $p$-parabolic equation 

\begin{subequations}\label{introeqex}
\begin{align}
 \frac{\partial u}{\partial t} - \nabla_x \cdot (|\nabla_x u|^{p-2}\nabla_x u)
&=f \quad
\mbox{in $\mathcal{D}'(Q_+)$}\\
u&=g \quad \mbox{on $ (\Omega \times \{0\}) \cup (\partial \Omega \times {\bf R}_+)$},
\end{align}
\end{subequations}
with $1<p<\infty$.

The boundary data is prescribed on the whole parabolic boundary, 
$ (\Omega \times \{0\}) \cup (\partial \Omega \times {\bf R}_+)$, and we study the problem 
of finding the ``largest possible'' classes of boundary and source data such that (\ref{introeq}) 
has a good meaning and is  uniquely solvable.

In the case of the elliptic  $p$-laplacian:
\begin{subequations}\label{introeqexellip}
\begin{align}
 - \nabla \cdot (|\nabla u|^{p-2}\nabla u)
 &=f \quad
\mbox{in $\mathcal{D}'(\Omega)$}\\
u&=g \quad \mbox{on $ \partial \Omega $},
\end{align}
\end{subequations}
it is well known that $W^{1,p}(\Omega)$ 
is a kind of golden mean. It has the useful property that:

Given $g \in W^{1,p}(\Omega)$, there
exists a unique solution $u \in W^{1,p}(\Omega)$ to the $p$-laplace
equation (\ref{introeqexellip}) such that  
$u-g$ belongs to the closure of $\mathcal D (\Omega)$ in the
$W^{1,p}(\Omega)$-norm topology.
Furthermore the source data ($f$ in (\ref{introeqexellip})) can then be taken as 
sums of first order derivatives of 
$L^{p/(p-1)}(\Omega)$-functions.

In this paper we construct an analogous optimal solution-space for equations
of the type (\ref{introeq}).

We point out that our results are new even in the linear case.
In the linear case, where $p=2$  and we denote $W^{s,2}$  by $H^s$, it is well known 
(see e.g. \cite{LN-MG} Vol. II) that the parabolic solution and lateral boundary value spaces, 
replacing the ``elliptic spaces'' $H^{s}(\Omega)$ and 
$H^{s-1/2}(\partial \Omega)$, are  $H^{s,s/2}(\Omega \times {\bf R}_+)$ and 
 $H^{s-1/2,s/2-1/4}(\partial \Omega \times {\bf R}_+)$. 
The initial data on $\Omega \times \{0\}$ should then belong to 
$H^{s-1}(\Omega)$ and the natural source data space is  $H^{s-2,s/2-1}(\Omega \times {\bf R}_+)$. 
With additional compatibility conditions for the 
coupling of the data in the ``corners'' of the space-time cylinder we then have  
unique solvability for the linear case when $s>1$ (see \cite{LN-MG}, Vol. II).
When $s=1$, the golden mean in the elliptic case,  several difficulties arise in 
the parabolic case.
One obvious difficulty is of course that we are in the borderline Sobolev 
imbedding case in the time direction (half-a-time derivative in 
$L^2({\bf R}_+,L^2(\Omega))$),
 and are thus for instance unable to define traces on 
$\Omega \times \{0\}$.

In  Theorem  \ref{th:mainlinear} we
give optimal results in the linear limiting case ($s=1$), and  a complete description 
of the space of solutions (compare 
with the non-optimal
 results in e.g. \cite{LN-MG},\cite{LD-SL-UR} and \cite{KP}). 

We use a  similar  construction of the solution space  
 (with new technical complications) in the non-linear case when $p \neq 2$.

Our solution space for a general $p$, $1<p<\infty$, (see Definition \ref{def:xspace}) is 
the sum of a Banach space carrying initial data and another Banach space 
carrying lateral boundary data.
It is  a dense subspace of the space of  $L^p(Q_+)$-functions, having
half order time derivatives in $L^2(Q_+)$ and first order space derivatives in
$L^p(Q_+)$.

This statement requires some explanation and the appropriate distribution theory, 
allowing fractional differentiation in the time direction of general $L^p$-functions in a 
space-time half cylinder, 
is developed. 
This analytic framework makes it possible to give a precise meaning to the 
fractional integration by parts for the time derivatives 
that is one of the key tools in our method.
We point out that we use  two different half-a-time derivatives (adjoint to each other) 
and that demanding these different derivatives  to belong to
$L^2(Q_+)$ 
 gives rise to  different function spaces.
 In Section 4 we investigate the relations between these different
 function spaces and discuss  some of their basic properties.
 It is for instance non-trivial to show that our function spaces are well behaved when we cut off 
(in a smooth way) in time. This is, apart from the fact that we are in the borderline 
Sobolev imbedding case in the time direction, due to the fact that they have non-homogeneous 
summability and regularity conditions, and that they are defined as spaces of distributions. 

Most of these technical problems arise already for functions defined
on the real line and half-line, and for clarity we have  moved
most of these arguments to an auxiliary section (Section 3) dealing
with this case.

The main result of this paper is Theorem \ref{th:mainnonlinear} which implies, among other things,
 that our solution space 
$X^{1,1/2}(Q_+)$
really is a true
analog of the space $W^{1,p}(\Omega)$ for the elliptic $p$-laplacian,
in the sense that:

Given $g \in X^{1,1/2}(Q_+)$ there
exists a unique solution $u \in X^{1,1/2}(Q_+)$ to the $p$-parabolic
equation (\ref{introeq})  such that  
$u-g$ belongs to the closure of $\mathcal D (Q_+)$ in the
$X^{1,1/2}(Q_+)$-norm topology.
Furthermore the source data ($f$ in (\ref{introeq})) can  be taken as 
sums of first order space derivatives of 
$L^{p/(p-1)}(Q_+)$-functions and half-a-time derivatives of $L^2(Q_+)$-functions.

For simplicity we shall assume throughout the paper 
that the boundary of $\Omega $ is smooth, but this assumption is 
only used to prove that we can regularize functions near the lateral boundary so that 
the different spaces of test functions we use are dense in the
corresponding function spaces (see Theorem \ref{th:densetestfunc}).

\section{Some analytical background.}
We will use the fractional calculus presented in \cite{F}. Here we first give 
a brief review of the notation and some results. We then extend the calculus to 
space-time half-cylinders in order to be able to discuss initial-boundary value 
problems.

The Fourier transform on the Schwartz class ${\mathcal S}({\bf R}^n, {\bf C})$ is defined by
\begin{equation}
\hat u (\xi) = \int_{{\bf R}^n} u(x)e^{-i2\pi x \cdot \xi }\, dx, \quad u \in 
{\mathcal S}({\bf R}^n, {\bf C}).
\end{equation}
The inverse will be denoted
\begin{equation}
\check u (\xi) = \int_{{\bf R}^n} u(x)e^{i2\pi x \cdot \xi} \, dx, \quad u \in 
{\mathcal S}({\bf R}^n, {\bf C}).
\end{equation}
The isotropic fractional Sobolev spaces are defined as follows.
\begin{definition}
For $s \in {\bf R}$ and $1<p<\infty$ let
\begin{equation}
H_p^s({\bf R}^n, {\bf C})= \{ u \in {\mathcal S'}({\bf R}^n, {\bf C});\;
((1+ |2\pi \xi |^2)^{s/2} 
 \hat u (\xi ))^{\vee}  \in L^p({\bf R}^n,{\bf C})  \}.
\end{equation}
\end{definition}
They are separable and reflexive Banach spaces 
with the obvious
norms. 
We will use the following  multi-index notation.
Let $\alpha =(\alpha_1, \dots ,\alpha_n) \in {\bf R}^n$ be an $n$-tuple.
We write $\alpha > 0$ if $\alpha_j > 0,\; j=1, \dots , n$;
$x^{\alpha}=x_1^{\alpha_1} \cdots x_n^{\alpha_n}$ when $x \in {\bf R}^n$;
$x^{\alpha}_+= {x_1^{\alpha_1}}_+\cdots  {x_n^{\alpha_n}}_+$,
 (where $t_+= \max (0,t)$ for $t \in {\bf R}$, with a similar definition
for $x^{\alpha}_-$) 
and $\Gamma (\alpha)= \Gamma(\alpha_1) \cdots \Gamma(\alpha_n)$, where
$\Gamma$ denotes the gamma function.
Furthermore we will sometimes write $k$ for the multi-index $(k,\dots,k)$,
the interpretation should be clear from the context.
We now define the classical Riemann-Liouville convolution operators.

\begin{definition}
For a multi-index $\alpha >0$, set
\begin{equation}
D_{\pm}^{-\alpha} u = \chi_{\pm}^{\alpha -1}* u,\quad u \in {\mathcal S}({\bf R}^n,{\bf C}),
\end{equation}
where the kernels  $\chi_{\pm}^{\alpha -1}$,
are given by
\begin{equation}
\chi_{\pm}^{\alpha -1} = \Gamma (\alpha)^{-1} (\cdot)_{\pm}^{\alpha-1}.
\end{equation}
\end{definition}
We extend the definition of $D_{\pm}^{\alpha}$ to general
multi-indices $\alpha \in {\bf R}^n$ in the usual way.
\begin{definition}
For $\alpha \in {\bf R}^n$ set
\begin{equation}
D_{\pm}^{\alpha}u= D^k D_{\pm}^{\alpha-k}u, \quad u 
\in {\mathcal S}({\bf R}^n,{\bf C}),
\end{equation}
where we choose the multi-index
$k \in \{0,1,2,\dots\}^n$ so that $k-\alpha >0$. 
\end{definition}
The definition is independent of the choice of $k$.

Although it is clear in this setting  how the support of a function is
affected under these mappings and also for instance that the operators 
map real valued
functions to real valued functions, other features become 
 transparent on the Fourier transform side.

Computing in  ${\mathcal S'}({\bf R}^n,{\bf C})$, we have 
for all
$\alpha \in {\bf R}^n $:
\begin{equation}
D_{\pm}^{\alpha}u = ((0{\pm}i 2 \pi \xi )^{\alpha} \hat u (\xi))^{\vee}, 
\quad u \in  {\mathcal S}({\bf R}^n,{\bf C}).
\end{equation}

We will use the following space of test functions.

\begin{definition}\label{def:testfunctions}
Let
\begin{multline}
\mathcal F ({\bf R}^n,{\bf C})  \\
=\left\{ u \in C^{\infty}({\bf R}^n,{\bf C});\quad 
\|u\|_{H^s_p({\bf R}^n,{\bf C})} < \infty,\;s \in {\bf R},\; 1< p<\infty 
\right\}.
\end{multline}
\end{definition}
$\mathcal F ({\bf R}^n,{\bf C})$ becomes a Fr\'echet space with the
topology generated by, for instance, the following  
family of semi-norms 
$\| \cdot \|_{H^s_p({\bf R}^n,{\bf C})},$ $  s \in \{0,1,2,\dots\},$ $
p=1+2^k,$ $ k\in {\bf Z} $.

We have the following dense continuous imbeddings,
\begin{equation}
\mathcal D ({\bf R}^n,{\bf C}) \hookrightarrow \mathcal S ({\bf R}^n,{\bf C}) 
\hookrightarrow \mathcal F ({\bf R}^n,{\bf C}) \hookrightarrow 
\mathcal E ({\bf R}^n,{\bf C}).
\end{equation}

An example of a function that belongs to $\mathcal F ({\bf R},{\bf C})$ but does not
belong to $\mathcal S ({\bf R},{\bf C})$ is $x \mapsto 1/(1+x^2)$.

For $\alpha \geq 0$ we now define the fractional derivatives
\begin{equation}
D_{\pm}^{\alpha} u = ((0 {\pm}  i 2 \pi \xi )^{\alpha} \hat u )^{\vee}, \qquad u \in 
\mathcal F ({\bf R}^n,{\bf C}).
\end{equation}
The operators $D_+^{\alpha}$ and $D_-^{\alpha}$ are adjoint to each other and 
they are connected through
the operator
\begin{equation}
H^{\alpha} = \prod_{k=1}^n (\cos(\pi \alpha_k) \mbox{Id} +
\sin(\pi \alpha_k) H_k),
\end{equation}
where Id is the identity operator and $H_k$ is the Hilbert transform
with respect to the $k$th variable, i.e.
\begin{equation}
H_k u(t)= \pi^{-1}\lim_{\epsilon \rightarrow +0} \int_{|s| \geq \epsilon}
\frac{u(t-se_k)}{s}\, ds, \quad u \in \mathcal F ({\bf R}^n,{\bf C}),
\end{equation}
where $e_k$ is the usual canonical $k$th basis vector in ${\bf R}^n$.
We have the following lemma.
\begin{lemma}\label{lemma:compositionrulesfunctions}
For $\alpha \geq 0$, $D_{\pm}^{\alpha}$ are continuous
linear operators on $\mathcal F ({\bf R}^n,{\bf C})$.
For $\alpha \in {\bf R}^n$, $H^{\alpha}$ is an isomorphism on
$\mathcal F ({\bf R}^n,{\bf C})$. For $\alpha, \beta \geq 0$ we have
\begin{eqnarray}
D_{\pm}^{\alpha}D_{\pm}^{\beta} = D_{\pm}^{\alpha+\beta},\\
D_+^{\alpha}H^{\alpha} =D_-^{\alpha}.
\end{eqnarray}
Furthermore all these operators commute on $\mathcal F ({\bf R}^n,{\bf C})$.
\end{lemma}
We note that for $\alpha \geq 0$
\begin{equation}\label{eq:adjointnessDD}
\int_{{\bf R}^n}  D_+^{\alpha} u \Phi  \, dx =\int_{{\bf R}^n} u D_-^{\alpha} \Phi \, dx,
\quad u,\Phi \in \mathcal F ({\bf R}^n,{\bf C}),
\end{equation}
and for $\alpha \in {\bf R}^n$
\begin{equation}\label{eq:adjointnessHH}
\int_{{\bf R}^n}  H^{\alpha} u \Phi \, dx =\int_{{\bf R}^n} u H^{-\alpha} \Phi \, dx,
\quad u,\Phi \in \mathcal F ({\bf R}^n,{\bf C}).
\end{equation}

Now let $\mathcal F' ({\bf R}^n,{\bf C})$ denote the space of continuous 
linear functionals
on $\mathcal F ({\bf R}^n,{\bf C})$, endowed with the weak$^*$ topology.

Inspired by (\ref{eq:adjointnessDD}) and (\ref{eq:adjointnessHH}), 
we extend the definition of 
$D_{\pm}^{\alpha}$ and $H^{\alpha}$ to $\mathcal F' ({\bf R}^n,{\bf C})$ by
duality in the obvious way.
\begin{definition}
For $u \in \mathcal F' ({\bf R}^n,{\bf C})$ and $\alpha \geq 0$ let
\begin{equation}
\langle D_{\pm}^{\alpha} u, \Phi \rangle 
:= \langle u , D_{\mp}^{\alpha} \Phi \rangle ,
\quad \Phi \in \mathcal F ({\bf R}^n,{\bf C}),
\end{equation}
and for $\alpha \in {\bf R}^n$ let
\begin{equation}
\langle H^{\alpha} u, \Phi \rangle := \langle u, H^{-\alpha} \Phi \rangle , 
\quad \Phi \in \mathcal F ({\bf R}^n,{\bf C}).
\end{equation}
\end{definition}
The counterpart of Lemma \ref{lemma:compositionrulesfunctions} is valid for
 $\mathcal F' ({\bf R}^n,{\bf C})$.
\begin{lemma}\label{lemma:compositionrulesdist}
For $\alpha \geq 0$, $D_{\pm}^{\alpha}$ are continuous
linear operators on $\mathcal F' ({\bf R}^n,{\bf C})$.
For $\alpha \in {\bf R}^n$, $H^{\alpha}$ is an isomorphism on
$\mathcal F' ({\bf R}^n,{\bf C})$. For $\alpha, \beta \geq 0$ we have
\begin{eqnarray}
D_{\pm}^{\alpha}D_{\pm}^{\beta} = D_{\pm}^{\alpha+\beta},\\
D_+^{\alpha}H^{\alpha} =D_-^{\alpha}.
\end{eqnarray}
Furthermore all these operators commute on $\mathcal F' ({\bf R}^n,{\bf C})$.
\end{lemma}

We recall that $D_{\pm}^{\alpha}$ and $H^{\alpha}$ all take real-valued 
functions
(distributions) to real-valued functions (distributions), and from now on all
functions and distributions will be real valued.
We will denote the subspaces of real-valued functions and distributions 
simply by
 $\mathcal F ({\bf R}^n)$ and  $\mathcal F' ({\bf R}^n)$.

In \cite{F} we studied parabolic operators on a space-time cylinder 
 $Q=\Omega \times {\bf R}$, where $\Omega$ was a connected and open set in 
${\bf R}^n$. We then introduced the following space of test functions.

\begin{definition}
Let 
$\mathcal F_{0,\cdot} (Q)$
denote the subspace of $\mathcal F ({\bf R}^n \times {\bf R})$ functions 
with support in 
$K \times {\bf R}$ for some compact subset $K \subset \Omega$.
\end{definition}

We put a pseudo-topology on $\mathcal F_{0,\cdot} (Q)$ by specifying what
sequential convergence means. We say that $\Phi_i \longrightarrow 0$ in
$\mathcal F_{0,\cdot} (Q)$ if and only if the supports of all $\Phi_i$'s are 
contained in a fixed set $K \times {\bf R}$, where $K \subset \Omega$ is a 
compact subset, and $\| D^{\alpha} \Phi_i \|_{L^P(Q)} \longrightarrow 0$ as
$i \longrightarrow \infty$ for all multi-indices $\alpha \in {\bf Z}_+^{n+1}$ and
$1<p<\infty$.

The corresponding space of distributions is then defined as follows.

\begin{definition}
If $u$ is a linear functional on $\mathcal F_{0,\cdot} (Q)$, then $u$ is in 
$\mathcal{F'}_{\cdot,\cdot} (Q)$ if and only if 
for every compact set $K \subset \Omega$,
there exist constants $C,p_1,\dots, p_N $ with $1<p_i<\infty,\quad i=1,\dots,N$
 and multi-indices $\alpha_1,\dots,\alpha_N$ with $\alpha_i \in {\bf Z}_+^{n+1}, \quad
i=1,\dots,N$ such that
\begin{equation}
|\langle u, \Phi \rangle| \leq C \sum_{i=1}^N \|D^{\alpha_i} \Phi \|_{L^{p_i}(Q)}
\end{equation}
for all $\Phi \in \mathcal F_{0,\cdot} (Q)$ with support in 
$K \times {\bf R}$.
\end{definition}

The motivation for these spaces is that they are invariant under  
fractional differentiation and Hilbert-transformation in the time 
variable, and ordinary differentiation in the space variables. In the
given topologies, these operations are continuous.

 For initial-boundary value problems, the parabolic operators 
will by defined on a space-time 
half-cylinder $Q_+=\Omega \times {\bf R}_+$, and
 we shall then need the following natural spaces of test functions  
defined on $Q_+$. 

{\bf Remark.} We shall use the same constructions on the real line and half-line, 
which can be thought of as the case $\Omega =\{0\}$ if we identify $\{0\} \times {\bf R}$ with 
${\bf R}$ and  $\{0\} \times {\bf R}_+$ with 
${\bf R}_+$.

\begin{definition}
Let 
$\mathcal F_{0,\cdot} (Q_+)$
denote the space of those functions defined on $Q_+$ that can 
be extended to all of $Q$ 
as elements in $\mathcal F_{0,\cdot} (Q)$.

Furthermore let 
$\mathcal F_{0,0} (Q_+)$
denote the space of those functions defined on $Q_+$ that can 
be extended by zero to all of $Q$ 
as elements in $\mathcal F_{0,\cdot} (Q)$.
\end{definition}
 (A zero in the first
position of course corresponds to zero boundary data on the lateral boundary and 
a
zero in the second position corresponds to zero initial data.)

By using the construction in \cite{SE} of a (total) extension operator, 
we see that $\mathcal F_{0,\cdot} (Q_+)$ can be identified with  
the space of all smooth functions $\Phi$, 
defined on $Q_+$, with support in $K \times R_+$  for some compact subset 
$K \subset \Omega$ (i.e.  they are zero on the complement, with respect to $Q_+$, 
of 
$K \times {\bf R}_+$), with 
$\| D^{\alpha} \Phi \|_{L^P(Q_+)} < \infty$ for all  
multi-indices $\alpha \in {\bf Z}_+^{n+1}$ and
$1<p<\infty$.

Thus, we can put an intrinsic pseudo-topology on $\mathcal F_{0,\cdot} (Q_+)$ by 
defining that  $\Phi_i \longrightarrow 0$ in
$\mathcal F_{0,\cdot} (Q_+)$ if and only if the supports of all $\Phi_i$ are 
contained in a fixed set $K \times {\bf R}_+$, where $K \subset \Omega$ is a 
compact subset, and $\| D^{\alpha} \Phi_i \|_{L^P(Q_+)} \longrightarrow 0$ as
$i \longrightarrow \infty$ for all multi-indices $\alpha \in {\bf Z}_+^{n+1}$ and
$1<p<\infty$.
Then $\mathcal F_{0,0} (Q_+)$ 
is a closed subspace of  $\mathcal F_{0,\cdot} (Q_+)$ with the induced topology.

We also note that $\mathcal D (Q_+)$ is densely continuously imbedded in 
$\mathcal F_{0,0} (Q_+)$.

Connected with these spaces of test functions are the following spaces of
distributions.

\begin{definition}
If $u$ is a linear functional on $\mathcal F_{0,\cdot} (Q_+)$, then $u$ is in 
$\mathcal{F'}_{\cdot,0} (Q_+)$ if and only if 
for every compact set $K \subset \Omega$,
there exist constants $C,p_1,\dots, p_N $ with $1<p_i<\infty,\quad i=1,\dots,N$
 and multi-indices $\alpha_1,\dots,\alpha_N$ with $\alpha_i \in {\bf Z}_+^{n+1}, \quad
i=1,\dots,N$ such that
\begin{equation}
|\langle u, \Phi \rangle| \leq C \sum_{i=1}^N \|D^{\alpha_i} \Phi \|_{L^{p_i}}
\end{equation}
for all $\Phi \in \mathcal F_{0,\cdot} (Q_+)$ with support in 
$K \times {\bf R}_+$.

Furthermore 
if $u$ is a linear functional on $\mathcal F_{0,0} (Q_+)$, then $u$ is in 
$\mathcal{F'}_{\cdot,\cdot} (Q_+)$ if and only if for every compact set 
$K \subset \Omega$,
there exist constants $C,p_1,\dots, p_N $ with $1<p_i<\infty,\quad i=1,\dots,N$
 and multi-indices $\alpha_1,\dots,\alpha_N$ with $\alpha_i \in {\bf Z}_+^{n+1}, \quad
i=1,\dots,N$ such that
\begin{equation}
|\langle u, \Phi \rangle| \leq C \sum_{i=1}^N \|D^{\alpha_i} \Phi \|_{L^{p_i}(Q_+)}
\end{equation}
for all $\Phi \in \mathcal F_{0,0} (Q_+)$ with support in 
$K \times {\bf R}_+$.
\end{definition}

The importance of these spaces comes from the fact that, 
 for a real-valued $\alpha \geq 0$, the operations
\begin{eqnarray}
\frac{\partial_+^{\alpha}}{\partial t^{\alpha}}:=D_+^{(0,\dots,0,\alpha)}: 
\mathcal{F}_{0,0}(Q_+) \longrightarrow 
\mathcal{F}_{0,0}(Q_+)\\
\frac{\partial_-^{\alpha}}{\partial t^{\alpha}}:=D_-^{(0,\dots,0,\alpha)}: 
\mathcal{F}_{0,\cdot}(Q_+) \longrightarrow 
\mathcal{F}_{0,\cdot}(Q_+)
\end{eqnarray}
are continuous. Ordinary differentiations with respect to the space variables are clearly
also continuous operations on these spaces.
We shall also use that the Hilbert-transform in the time variable
\begin{equation}
h:= H^{(0,\dots,0,1/2)}: 
\mathcal{F}_{0,0}(Q_+) \longrightarrow 
\mathcal{F}_{0,\cdot}(Q_+),
\end{equation}
is a continuous operator.

Extending these operators by duality in the obvious way we get that
\begin{eqnarray}
\frac{\partial_+^{\alpha}}{\partial t^{\alpha}}: \mathcal{F'}_{\cdot,0}(Q_+) \longrightarrow 
\mathcal{F'}_{\cdot,0}(Q_+),\\
\frac{\partial_-^{\alpha}}{\partial t^{\alpha}}: \mathcal{F'}_{\cdot,\cdot}(Q_+) \longrightarrow 
\mathcal{F'}_{\cdot,\cdot}(Q_+),\\
h:
\mathcal{F'}_{\cdot,0}(Q_+) \longrightarrow 
\mathcal{F'}_{\cdot,\cdot}(Q_+),
\end{eqnarray}
and taking ordinary derivatives in the space variables, are continuous operations.

Using the total extension operator from \cite{SE},  one can show that we can 
identify  $\mathcal{F'}_{\cdot,0}(Q_+)$ with the 
space of $\mathcal{F'}_{\cdot,\cdot}(Q)$-distributions that are zero on 
$\Omega \times (-\infty,0)$. 

Since $\mathcal D (Q_+)$ is densely continuously imbedded in 
$\mathcal{F}_{0,0}(Q_+)$, we get  
that $\mathcal{F'}_{\cdot,\cdot}(Q_+)$ is a continuously imbedded 
subspace of $\mathcal {D'}(Q_+)$.

We remark that the space $\mathcal{F'}_{\cdot,0}(Q_+)$ contains elements supported on 
$\Omega \times \{ 0 \}$. In fact
\begin{equation}
\mathcal{F'}_{\cdot,\cdot}(Q_+) \simeq \mathcal{F'}_{\cdot,0}(Q_+)/
\mathcal{F^{\circ}}_{0,0}(Q_+),
\end{equation}
where $\mathcal{F^{\circ}}_{0,0}(Q_+)= 
\left\{ \xi \in \mathcal{F'}_{\cdot,0}(Q_+);\quad
\langle \xi,\Phi \rangle =0,\; \Phi \in \mathcal{F}_{0,0}(Q_+) \right\} $. 

Finally, since $\mathcal{F}_{0,0}(Q_+)$ is densely continuously imbedded in 
$L^p(Q_+)$ when $1<p<\infty$, clearly 
$L^p(Q_+)$ is continuously imbedded in  both
$\mathcal{F'}_{\cdot,\cdot}(Q_+)$ and $\mathcal{F'}_{\cdot,0}(Q_+)$
when $1<p<\infty$. Thus
\begin{eqnarray}
\frac{\partial_+^{\alpha}}{\partial t^{\alpha}}: L^p(Q_+) \longrightarrow 
\mathcal{F'}_{\cdot,0}(Q_+)\\
\frac{\partial_-^{\alpha}}{\partial t^{\alpha}}: L^p(Q_+)\longrightarrow 
\mathcal{F'}_{\cdot,\cdot}(Q_+),
\end{eqnarray}
are well-defined continuous operations when $1<p<\infty$.

\section{Auxiliary spaces on the real line and half-line.}

We shall use the following auxiliary spaces defined on ${\bf R}$ and in the 
definition $\frac{\partial^{1/2}_- }{\partial t^{1/2}}$ should be understood
in the $\mathcal F' ({\bf R})$ distribution sense.

\begin{definition}
For $1<p<\infty$, set
\begin{equation}
B^{1/2}({\bf R})=\left\{ u \in L^p({\bf R});\; 
\frac{\partial^{1/2}_- u}{\partial t^{1/2}} 
\in L^2({\bf R})\right\}.
\end{equation}
\end {definition} 

We equip these spaces with the following norms.
 
\begin{equation}
\|u\|_{B^{1,1/2}({\bf R})}:= 
\|\frac{\partial^{1/2}_- u}{\partial t^{1/2}} \|_{L^2({\bf R})} + 
\|u\|_{L^p({\bf R})}.
\end{equation}

Computing in $\mathcal {F'}({\bf R})$ we see that
we can represent these spaces as closed subspaces of 
the direct sums 
$L^2({\bf R})\oplus L^p({\bf R})$, and thus they are 
reflexive and separable Banach spaces in the topologies 
arising from the given norms.

If $\{\psi_{\epsilon}\}$ is a regularizing  sequence it is clear that
\begin{equation}
\|\psi_{\epsilon} \ast u \|_{B^{1/2}({\bf R})} \leq \|u\|_{B^{1/2}({\bf R})} \;, 
\end{equation}
and thus smooth functions are dense in $B^{1/2}({\bf R})$.
 
Due to the definition using distributions and to the inhomogeniety of our 
summability conditions, it is unfortunately 
not so easy to cut off 
in time and in this way show that $\mathcal{F}({\bf R})$ (or $\mathcal{D}({\bf R})$) 
is dense in  $B^{1/2}({\bf R})$. 
Nevertheless this is true.

\begin{lemma}
 The space of testfunctions $\mathcal{F}({\bf R})$ is dense in 
 $B^{1/2}({\bf R})$.
\end{lemma}
{\bf Proof.}
The proof is based on a non-linear version of the Riesz representation theorem.

We (temporarily) denote the closure of $\mathcal{F}({\bf R})$ in  $B^{1/2}({\bf R})$
by  $B^{1/2}_0({\bf R})$, and we shall show that 
 $B^{1/2}_0({\bf R}) = B^{1/2}({\bf R})$.

Set
\begin{equation}
T(u)= \frac{\partial u}{\partial t} +|u|^{p-2}u.
\end{equation}

By fractional integration by parts 
\begin{equation}
\langle T(u), \Phi \rangle = \int_{{\bf R}} 
\frac{\partial_+^{1/2}u}{\partial t^{1/2}}
\frac{\partial_-^{1/2}\Phi}{\partial t^{1/2}} +|u|^{p-2}u \Phi \, dt \; ; \quad 
\Phi \in \mathcal F({\bf R}),
\end{equation}
 and H\"older's inequality, it is clear that 
\begin{equation}
T :  B^{1/2}({\bf R}) \longrightarrow B^{1/2}_0({\bf R})^*.
\end{equation}
is continuous.

We notice that   
\begin{equation}
T :  B^{1/2}_0({\bf R}) \longrightarrow B^{1/2}_0({\bf R})^*,
\end{equation}
is weakly continuous and monotone (for definitions see [KS] or \cite{F}).

By M. Riesz' conjugate function theorem, which says that 
the Hilbert transform 
$h$ is bounded from $L^p({\bf R})$ to  $L^p({\bf R})$ (recall that $1<p<\infty$), 
we see that
the operators $H^{\alpha}$ introduced above are isomorphisms on 
 $B^{1/2}_0({\bf R})$.

Now for any $\alpha \in (0,1/2)$ we have 
\begin{eqnarray}
\langle T(u), H^{-\alpha}(u) \rangle \geq 
 \int_{{\bf R}} 
\sin (\pi \alpha) \frac{\partial_+^{1/2}u}{\partial t^{1/2}}
\frac{\partial_+^{1/2}u}{\partial t^{1/2}}\\
 +(\cos (\pi \alpha) 
-\sin (\pi \alpha) C)|u|^p \, dt \quad ; u \in \mathcal F ({\bf R}),
\end{eqnarray}
where $C <\infty$ is a constant such that
\begin{equation}
\| h(u)\|_{L^p({\bf R})} \leq C \| u\|_{L^p({\bf R})}.
\end{equation}
Choosing $\alpha \in (0,1/2)$ small enough we see that $H^{\alpha} \circ T$ 
is coercive. 
It follows that 
$T$ is 
a bijection (see \cite{F}  for  
this functional-analytic result and similar arguments).

Thus given $u \in  B^{1/2}({\bf R})$ there exists a unique 
 $v \in  B^{1/2}_0({\bf R})$ such that $T(u)=T(v)$ in $\mathcal F'({\bf R})$, i.e.
\begin{equation}\label{eq:modelopeq}
\frac{\partial (u-v)}{\partial t} +(|u|^{p-2}u -|v|^{p-2}v) =0.
\end{equation}

This shows that the difference of elements with the same image has more 
regularity in time, namely $\frac{\partial (u-v)}{\partial t} \in L^{p/(p-1)}({\bf R})$.

The class of $L^p({\bf R})$ functions with derivatives in  $L^{p/(p-1)}({\bf R})$
is stable under regularization and thus by 
a continuity argument we see that we can test with $\chi(u-v)$,
 where $\chi$ is a cut off 
function in time, in equation (\ref{eq:modelopeq}).  We get that 
(for a canonical continuous representative)
$t \mapsto |u-v|(t)$ is decreasing. Since $u-v$ belongs to $L^p({\bf R})$, 
 we conclude that $u=v$. 
The lemma follows.
$\Box$

We are now in position to prove the following lemma.

\begin{lemma}\label{lem:eqnorm}
If  $u \in  B^{1/2}({\bf R})$ then
\begin{equation}
\iint_{{\bf R}\times{\bf R}} \left|\frac{u(s)-u(t)}{s-t}\right|^2 \, ds\,dt 
=2 \pi \int_{{\bf R}} \left|\frac{\partial_-^{1/2} u}{\partial t^{1/2}} \right|^2 \, dt.
\end{equation}
\end{lemma}

{\bf Proof.}
Since  $\mathcal{F}({\bf R})$ is dense in 
 $B^{1/2}({\bf R})$ we can compute using the Fourier transform.
\begin{eqnarray}
 \int_{{\bf R}} \left|\frac{\partial_-^{1/2} u}{\partial t^{1/2}} \right|^2 \, dt
=\int_{{\bf R}} 2\pi |\tau||\hat{u}|^2 \,d\tau \\
=\frac{1}{2 \pi} \iint_{{\bf R}\times {\bf R}}\frac{|1- e^{i2\pi \tau s}|^2}{s^2}
|\hat{u}(\tau)|^2 \,d\tau\,ds.
\end{eqnarray}
Using Parseval's formula the lemma follows.
$\Box$

We note the following scaling and translation invariance 
\begin{eqnarray}\label{eq:scaletransinv}\nonumber
\iint_{{\bf R}\times{\bf R}} 
\left|\frac{u(a(s-b))-u(a(t-b))}{s-t}\right|^2 \, ds\,dt\\
=\iint_{{\bf R}\times{\bf R}} \left|\frac{u(s)-u(t)}{s-t}\right|^2 
\, ds\,dt \; ; a,b \in {\bf R}.
\end{eqnarray} 

We also note the following fact.
\begin{lemma}
The space  $B^{1/2}({\bf R})$ is continuously imbedded in the space of 
functions with vanishing mean oscillation, $VMO({\bf R})$.
\end{lemma}
{\bf Proof.}
Let $I \subset {\bf R}$ denote a bounded  interval and let $u_I$ denote the 
mean value of $u \in B^{1/2}({\bf R})$ over $I$. Then by Jensen's 
inequality 
\begin{equation}\label{eq:vmoh1/2}
\frac{1}{|I|} \int_I |u-u_I|^2\,dt 
\leq \iint_{I \times I} \left|\frac{u(s)-u(t)}{s-t}\right|^2\,ds\,dt.
\end{equation}
$\Box$

Using the form of the norm in Lemma \ref{lem:eqnorm}, we can now show that we have good 
estimates in  the $B^{1/2}({\bf R})$-norm for the following cut-off operation.

\begin{lemma}\label{lem:cutoff}
 Let $\chi_n$ be the piecewise affine function 
that is one on $(-n,n)$, zero on  $(-\infty,-2n) \cup (2n,\infty)$ and affine in between. 
Let $I_n=(-2n,2n)$ and for 
 $u \in B^{1/2}({\bf R})$, denote the mean value of $u$ over $I_n$ by $u_{I_n}$. Then 
there exists a constant $C$ such that
\begin{eqnarray}\label{eq:cutoff}\nonumber
\iint_{{\bf R}\times{\bf R}} \left|\frac{\chi_n(u-u_{I_n})(s)-\chi_n(u-u_{I})(t)}{s-t}\right|^2 
\, ds\,dt \\
\leq C \iint_{{\bf R}\times{\bf R}} \left|\frac{u(s)-u(t)}{s-t}\right|^2 
\, ds\,dt \; ,\\
\| \chi_n(u-u_{I_n})\|^p_{L^p({\bf R})} \leq  C \| u \|^p_{L^p({\bf R})}
\quad ; u \in B^{1/2}({\bf R}).
\end{eqnarray}
Furthermore 
$\chi_n(u-u_{I_n}) \rightarrow u$ in $B^{1/2}({\bf R})$ as $n \longrightarrow \infty$.
\end{lemma}

{\bf Proof.}
The boundedness of the cut-off operation in the $L^p$-norm follows from Jensen's inequality.
For the $L^2$-part of the norm an elementary computation gives us
\begin{eqnarray}\nonumber
\iint_{{\bf R}\times{\bf R}} \left|\frac{\chi_n(u-u_{I_n})(s)-\chi_n(u-u_{I})(t)}{s-t}\right|^2 
\, ds\,dt \\
\leq C \left\{ \frac{1}{|I_n|} \int_{I_n} |u-u_{I_n}|^2\,dt +
\iint_{{\bf R}\times{\bf R}} \left|\frac{u(s)-u(t)}{s-t}\right|^2 
\, ds\,dt \right\},
\end{eqnarray}
and thus (\ref{eq:cutoff}) follows using (\ref{eq:vmoh1/2}).
That $\chi_n u \rightarrow u$ in $L^p({\bf R})$ is clear.
If $u$ has compact support, 
since $p>1$, using Jensen's inequality, 
we see that $\chi_n u_{I_n} \rightarrow 0$ in $L^p({\bf R})$.
Since by Jensen's inequality   $\chi_n u_{I_n}$ is uniformly bounded in $L^p({\bf R})$, a density 
argument proves that $\chi_n(u-u_{I_n}) \rightarrow u$ in $L^p({\bf R})$.
That $\chi_n(u-u_{I_n}) \rightarrow u$ for the $L^2$-part of the norm  
follows since by an elementary computation
\begin{eqnarray}
\iint_{{\bf R}\times{\bf R}} \left|\frac{(1-\chi_n)(u-u_{I_n})(s)-(1-\chi_n)(u-u_{I})(t)}{s-t}\right|^2 
\, ds\,dt \\
\leq C \left\{ \frac{1}{|I_n|} \int_{I_n} |u-u_{I_n}|^2\,dt +
\iint_{|t|>n} \left|\frac{u(s)-u(t)}{s-t}\right|^2 
\, ds\,dt \right\}.
\end{eqnarray}
The last term clearly tends to zero as $n$ tends to infinity. We only have to prove that also
\begin{equation}
\frac{1}{|I_n|} \int_{I_n} |u-u_{I_n}|^2\,dt \longrightarrow 0
\end{equation}
as $n \rightarrow \infty$.
This is true since 
\begin{eqnarray}\nonumber
\frac{1}{|I_n|} \int_{I_n} |u-u_{I_n}|^2\,dt 
\leq \frac{1}{4n^2}\iint_{I_n \times I_n} |u(s)-u(t)|^2\,ds\,dt\\\nonumber
\leq C \left\{ \frac{\log^2 n}{n^2} \iint_{|s|,|t| \leq \log n}   
\left|\frac{u(s)-u(t)}{s-t}\right|^2 \, ds \, dt \right.\\ 
 \left.+  \iint_{|t| \geq \log n}   
\left|\frac{u(s)-u(t)}{s-t}\right|^2  \, ds \, dt \right\},
\end{eqnarray}
which clearly tends to zero as $n$ tends to infinity.
$\Box$

{\bf Remark.}
We subtracted the mean value in the argument above in order not to have to
rely on the fact that $u \in L^p({\bf R})$ when proving boundedness for the half-derivatives.
This is crucial when we later use the same argument  on  functions defined in a space-time cylinder.
In preparation for this we also note that, by regularizing, the lemma gives us an explicit sequence of 
$\mathcal{D} ({\bf R})$-functions tending to a given element in  $B^{1/2}({\bf R})$.

We now introduce two sets of spaces defined on the real half-line.

\begin{definition}
Let 
$ B^{1/2}_0({\bf R}_+)$ be the space of functions defined on ${\bf R}_+$ 
that can be extended by zero as elements in  $B^{1/2}({\bf R})$.

Furthermore
let 
$ B^{1/2}({\bf R}_+)$ be the space of functions defined on ${\bf R}_+$ 
that can be extended as elements in  $B^{1/2}({\bf R})$.
\end{definition}

{\bf Remark.}
The space $ B^{1/2}_0({\bf R}_+)$ can of course be identified 
with the closed subspace of $ B^{1/2}({\bf R})$ of functions with support in ${\bf R}_+$.

We now give two simple lemmas, giving intrinsic descriptions of $
B^{1/2}_0({\bf R}_+)$
 and $ B^{1/2}({\bf R}_+)$. We omit the proofs, which are
 straightforward elementary  computations 
 using the form of the norm
 in Lemma \ref{lem:eqnorm}.

\begin{lemma}\label{lem:bhalf0norm}
  The function space $B^{1/2}_0({\bf R}_+)$ is precisely the set of 
$L^p({\bf R}_+)$-functions such that the following norm is bounded:
\begin{eqnarray}\nonumber\label{eq:bhalf0norm}
\| u\|_{B^{1/2}_0({\bf R}_+)}:=\|u\|_{L^p({\bf R}_+)} +
\left\{ \int_{{\bf R}_+} \frac{u^2(t)}{t} \,dt \right.\\
\left. +
\iint_{{\bf R}_+ \times {\bf R}_+} 
\left(\frac{u(s)-u(t)}{s-t}\right)^2 \, ds\,dt \right\}^{1/2}.
\end{eqnarray}
\end{lemma}

\begin{lemma}
  The function space $B^{1/2}({\bf R}_+)$ is precisely the set of 
$L^p({\bf R}_+)$-functions such that the following norm is bounded:
\begin{equation}\label{eq:bhalfnorm}
\| u\|_{B^{1/2}({\bf R}_+)}:=\|u\|_{L^p({\bf R}_+)} +\left\{
\iint_{{\bf R}_+ \times {\bf R}_+} 
 \left(\frac{u(s)-u(t)}{s-t}\right)^2\, ds\,dt \right\}^{1/2}.
\end{equation}
Furthermore,  a continuous  symmetric extension operator  
from  $B^{1/2}({\bf R}_+)$ to  $B^{1/2}({\bf R})$
is given by  $E_S(u)(t)=u(|t|)$.
\end{lemma}

We have the following density results:
\begin{lemma}\label{lem:densereal}
The space $\mathcal F ({\bf R}_+)$ is dense in $B^{1/2}({\bf R}_+)$ and 
$\mathcal F_0 ({\bf R}_+)$ is dense in $B^{1/2}_0({\bf R}_+)$.
\end{lemma}
{ \bf Proof.}
That  $\mathcal F ({\bf R}_+)$ is dense in $B^{1/2}({\bf R}_+)$ follows immidiately 
from the fact that  $\mathcal F ({\bf R})$ is dense in $B^{1/2}({\bf R})$.
The argument to prove that $\mathcal F_0 ({\bf R}_+)$ is dense in $B^{1/2}_0({\bf R}_+)$ 
is a little more delicate. Given  $u \in B^{1/2}_0({\bf R}_+)$, apriori we only know 
that there exists a sequence of testfunctions in  $\mathcal F ({\bf R})$ approaching 
$u$ in  the $B^{1/2}({\bf R})$-norm.

Given $u \in B^{1/2}_0({\bf R}_+)$ we will show that we can cut-off.  
Let $\chi_n$ be the piecewise affine function that is one on $(0,n)$, zero on 
$(2n,\infty)$ and affine in between. We will show that $\chi_n u \longrightarrow u$ 
in  $B^{1/2}_0({\bf R}_+)$. Taking this for granted we can regularize with a 
regularizing sequence having support in ${\bf R}_+$   which gives us the lemma.

That $\chi_n u \longrightarrow u$ in $L^p({\bf R}_+)$ is clear.
We now estimate the $L^2$-part of the norm.  An elementary computation gives us
\begin{eqnarray}\nonumber
 \int_{{\bf R}_+} \frac{((1-\chi_n)u)^2(t)}{t} \,dt 
+ \iint_{{\bf R}_+ \times {\bf R}_+} 
\frac{((1-\chi_n)u(s)-(1-\chi_n)u(t))^2}{(s-t)^2} \, ds\,dt \\ \nonumber
\leq C \left\{ \frac{1}{2n} \int_0^{2n} u^2(t) \,dt +
\iint_{(n,\infty) \times {\bf R}_+} 
\left(\frac{u(s)-u(t)}{s-t}\right)^2\, ds\,dt \right.\\
\left. +\int_n^{\infty} \frac{u^2(t)}{t} \,dt   \right\}.\quad
\end{eqnarray}
The last two terms above clearly tend to zero as $n \rightarrow \infty$. 
To estimate the first term, we integrate by parts
(we may assume that $u$ is smooth, it is the decay at 
infinity that is the issue).
\begin{eqnarray}\nonumber
 \frac{1}{2n} \int_0^{2n} u^2(t) \,dt 
= \frac{1}{2n} \int_0^{2n} \left( \int_0^{2n} \frac{u^2(s)}{s}\, ds -
 \int_0^{t} \frac{u^2(s)}{s}\, ds \right) \, dt\\
\leq   \frac{1}{2n} \int_{\log n}^{2n}  \int_t^{2n}
  \frac{u^2(s)}{s} \,ds  \, dt 
+ \frac{\log n}{2n} \int_0^{2n} \frac{u^2(s)}{s} \, ds,
\end{eqnarray}
which clearly tends to zero as $n$ tends to infinity. The lemma follows.
$\Box$

We now give the following equivalent
characterization of  $B^{1/2}_0({\bf R}_+)$.

\begin{lemma}
A function $u \in L^p({\bf R}_+)$ belongs to $B^{1/2}_0({\bf R}_+)$ if and only if
the $\mathcal F'_0 ({\bf R}_+)$-distribution derivative 
$\frac{\partial^{1/2}_+ u}{\partial t^{1/2}} 
\in L^2({\bf R}_+)$. Furthermore an equivalent norm on $B^{1/2}_0({\bf R}_+)$ 
is given by
\begin{equation}
\|u\| = \| u\|_{L^p({\bf R}_+)} + 
\|\frac{\partial^{1/2}_+ u}{\partial t^{1/2}}\|_{L^2({\bf R}_+)}.
\end{equation}
\end{lemma}
{\bf Remark.} We recall that the 
$\mathcal F'_0 ({\bf R}_+)$-distribution derivative, apart from what happens inside 
${\bf R}_+$, also controls what happens on the boundary $\{0\}$. 
The fact that $\frac{\partial^{1/2}_+ u}{\partial t^{1/2}} 
\in L^2({\bf R}_+)$ thus actually contains a lot of information about 
$u$'s behaviour at $0$.

{\bf Proof.}

It is clear that a function in  $B^{1/2}_0({\bf R}_+)$ has the
$\mathcal F'_0 ({\bf R}_+)$-distribution derivative 
$\frac{\partial^{1/2}_+ u}{\partial t^{1/2}}$ 
in $L^2({\bf R}_+)$. 

On the other hand, 
let $E_0$ be the extension by zero operator. Then  if $u \in L^p({\bf
  R}_+)$ and 
the $\mathcal F'_0 ({\bf R}_+)$-distribution derivative 
$\frac{\partial^{1/2}_+ u}{\partial t^{1/2}} 
\in L^2({\bf R}_+)$ we have
\begin{equation}
\int_{\bf R} E_0(u) \frac{\partial^{1/2}_-\Phi}{\partial t^{1/2}}\,dt
=\int_{\bf R}  E_0(\frac{\partial^{1/2}_+ u}{\partial t^{1/2}}) \Phi \,dt\;
;\; \Phi \in \mathcal F ({\bf R}).
\end{equation}
This shows that the  $\mathcal F' ({\bf R})$-distribution derivative 
$\frac{\partial^{1/2}_+ E_0(u)}{\partial t^{1/2}} $ belongs to
$L^2({\bf R})$.
An easy computation shows that
\begin{eqnarray}\nonumber
   \int_{{\bf R}_+} 
\left| \frac{\partial_+^{1/2} u}{\partial t^{1/2}} \right|^2 \, dt =
 \int_{{\bf R}} 
\left| \frac{\partial_+^{1/2} E_0(u)(t)}{\partial t^{1/2}} \right|^2 \,
dt \\
\sim
\iint_{ {\bf R}_+ \times {\bf R}_+ } 
\left| \frac{E_0(u)(s)-E_0(u)(t)}{s-t}\right|^2 \, ds\,dt
+\int_{\bf R_+} \frac{E_0(u)^2}{t}\, dt.
  \end{eqnarray}
  Since $E_0(u)=u$ on $(0,\infty)$ the lemma follows.
  $\Box$

We now give a corresponding equivalent norm on 
$B^{1/2}({\bf R}_+)$.

\begin{lemma}
If $u \in B^{1/2}({\bf R}_+)$, 
then the $\mathcal F' ({\bf R}_+)$-distribution derivative 
$\frac{\partial^{1/2}_- u}{\partial t^{1/2}}$ 
belongs to $L^2({\bf R}_+)$. 

Furthermore an equivalent norm on $B^{1/2}({\bf R}_+)$ 
is given by
\begin{equation}
\|u\| = \| u\|_{L^p({\bf R})} + 
\|\frac{\partial^{1/2}_- u}{\partial t^{1/2}}\|_{L^2({\bf R}_+)}.
\end{equation}
\end{lemma}
{\bf Remark.} In contrast to the  $\mathcal F'_0 ({\bf R}_+)$-distribution 
derivative, the  $\mathcal F' ({\bf R}_+)$-distribution derivative that we use 
in this definition ``does not see'' what happens on the boundary, $\{0\}$. 

{\bf Proof.}
 Since $\mathcal F ({\bf R}_+)$ is dense in  $B^{1/2}({\bf R}_+)$, it is enough to show that 
\begin{equation}
 \int_{{\bf R}_+} 
\left| \frac{\partial_-^{1/2} u}{\partial t^{1/2}} \right|^2 \, dt \sim
\iint_{ {\bf R}_+ \times {\bf R}_+ } 
\left| \frac{u(s)-u(t)}{s-t}\right|^2 \, ds\,dt ,
\end{equation}
for functions in  $\mathcal F ({\bf R}_+)$, where $\sim$ means that 
the seminorms are equivalent.

For $p=2$ we (temporarily) denote the closure of  
$\mathcal F ({\bf R}_+)$ in the norm
\begin{equation}
\|u \| =\|\frac{\partial_-^{1/2} u}{\partial t^{1/2}} \|_{L^2({\bf R}_+)} 
+\|u\|_{L^2({\bf R}_+)},
 \end{equation}
 by $H$.

It follows directely from the definitions, and the fact that 
 $\mathcal F ({\bf R}_+)$ is dense in $B^{1/2}({\bf R}_+)$, that 
$B^{1/2}({\bf R}_+)$ is continuously imbedded in $H$.

 We shall now show that in fact $H=B^{1/2}({\bf R}_+)$. 

Let $T$ denote the operator 
$T:u \mapsto \frac{\partial u}{\partial t} +u$. 
Then  
$
T: B^{1/2}_0({\bf R}_+) \longrightarrow H^*
$
is continuous. This follows from fractional integration by parts,
\begin{eqnarray}\nonumber
\langle Tu, \Phi \rangle = \left( \frac{\partial_+^{1/2} u}{\partial t^{1/2}},
\frac{\partial_-^{1/2} \Phi}{\partial t^{1/2}} \right)_{L^2} 
+\left(u,\Phi\right)_{L^2} \\
;\quad \Phi \in \mathcal F ({\bf R}_+), \;  u \in \mathcal F_0 ({\bf R}_+),
\end{eqnarray}
and the fact that  $\mathcal F ({\bf R}_+)$ is dense in  $H$ and  that 
$\mathcal F_0 ({\bf R}_+)$ is dense in 
$B^{1/2}_0({\bf R}_+)$.

Now by the Hahn-Banach theorem,  given $\xi \in H^*$ there exist
elements $u, v \in L^2({\bf R}_+)$ such that
\begin{equation}
\langle \xi, \Phi \rangle = 
\left(u,\frac{\partial_-^{1/2} \Phi}{\partial t^{1/2}} \right)_{L^2} +
\left(v,\Phi\right)_{L^2}
;\quad \Phi \in \mathcal F ({\bf R}_+).
\end{equation}

We can thus extend $\xi$ by zero to an element $E_0(\xi)$ of $B^{1/2}({\bf R})^*$.
Since $ T:B^{1/2}({\bf R}) \rightarrow B^{1/2}({\bf R})^*$ is an isomorphism, 
we can find a unique element   $ u \in B^{1/2}({\bf R})$ such that
$Tu = E_0(\xi)$ in $\mathcal {F'} ({\bf R})$. But this holds if and only if 
 $ u \in B^{1/2}_0({\bf R}_+)$ and 
$Tu = \xi$ in $\mathcal F'_0 ({\bf R}_+)$.
 
Thus $T: B^{1/2}_0({\bf R}_+) \longrightarrow H^*$
is an isomorphism.

Furthermore, by direct computation (or by interpolation (recall that $p=2$)), we know that
\begin{equation}
T: B^{1/2}_0({\bf R}_+) \longrightarrow  B^{1/2}({\bf R}_+)^*
\end{equation}
is an isomorphism.

Since $\mathcal {F} ({\bf R}_+)$ is densely continuously imbedded in 
both  $H$ and $B^{1/2}({\bf R}_+)$ and thus   
$H^*$ and $B^{1/2}({\bf R}_+)^*$ both are well defined subspaces in 
 $\mathcal F'_0 ({\bf R}_+)$, we see that 
 $H^*$ and $B^{1/2}({\bf R}_+)^*$ are identical as 
subspaces of  $\mathcal F'_0 ({\bf R}_+)$ and equivalent as Hilbert spaces.

Since $B^{1/2}({\bf R}_+) \hookrightarrow H$, by Riesz representation theorem, 
this implies that  $H$ and $B^{1/2}({\bf R}_+)$ have  equivalent norms.

From a scaling argument it now follows that
\begin{equation}
 \int_{{\bf R}_+} 
\left| \frac{\partial_-^{1/2} u}{\partial t^{1/2}} \right|^2 \, dt \sim
\iint_{ {\bf R}_+ \times {\bf R}_+ } 
\left| \frac{u(s)-u(t)}{s-t}\right|^2 \, ds\,dt ,
\end{equation}
for functions in  $B^{1/2}({\bf R}_+)$.
The lemma follows.
$\Box$

\section{Parabolic Equations.}
We shall consider operators of the form
\begin{equation}\label{eq:T}
Tu = \frac{\partial u}{\partial t} - \nabla_x\cdot  A(x,t,\nabla_x u),
\end{equation}
on a space-time cylinder $Q_+=\Omega \times {\bf R_+}$, where $\Omega$ is
an open and bounded set in ${\bf R}^n$ with smooth boundary.

We shall assume the following structural conditions for
the function 
$A: \Omega \times {\bf R}_+ \times {\bf R}^n \longrightarrow
{\bf R}^n$.
\begin{enumerate}
\item $Q_+ \ni (x,t) \mapsto A(x,t,\xi)$ is Lebesgue measurable
for every fixed $\xi \in {\bf R}^n$.
\item ${\bf R}^n \ni \xi \mapsto A(x,t,\xi)$ is continuous
for almost every $(x,t) \in Q_+$.
\item For every $\xi,\eta  \in {\bf R}^n, \xi \neq \eta$ and almost every
$(x,t) \in Q_+$, we have
\begin{equation}\label{eq:strictmonotonicity}
(A(x,t,\xi)-A(x,t,\eta),\xi -\eta) > 0.
\end{equation}
\item There exists $p \in (1,\infty)$, a constant $\lambda >0$
and a function $ h \in L^1(Q_+)$ such that for every $\xi \in {\bf R}^n$
and almost every$(x,t) \in Q_+$:
\begin{equation}
(A(x,t,\xi), \xi) \geq \lambda |\xi|^p - h(x,t).
\end{equation}
\item There exists a constant $\Lambda \geq \lambda >0$ and a
function $H \in L^{p/(p-1)}(Q_+)$ such that for every $\xi \in {\bf R}^n$
and almost every $(x,t)\in Q_+$:
\begin{equation}\label{eq:boundednesscond}
|A(x,t,\xi)| \leq \Lambda |\xi|^{p-1} + H(x,t).
\end{equation}
\end{enumerate}
The Carath\'eodory conditions 1 and 2 above guarantee that the function
$Q \ni (x,t) \mapsto A(x,t,\Phi(x,t))$ is measurable for every function
$\Phi \in L^p(Q_+,{\bf R}^n)$.
Condition 3 is a strict monotonicity condition that gives us
uniqueness results.
Conditions 4 (coercivity) and 5 (boundedness) give us apriori
 estimates that imply
existence results (see \cite{F}).

We now introduce some function spaces, and in their definitions 
$\partial_-^{1/2} /\partial t^{1/2}$ should be understood in the 
$\mathcal{F'}_{\cdot,\cdot}(Q)$ distribution-sense.

\begin{definition}\label{def:functionspaces}
For $1<p<\infty$, set
\begin{eqnarray}\nonumber
B^{1,1/2}_{\cdot,\cdot}(Q)=\left\{ u \in L^p(Q);\; 
\frac{\partial^{1/2}_- u}{\partial t^{1/2}} 
\in L^2(Q)\right. \\
, \left. \frac{\partial u}{\partial x_i} \in L^p(Q),
\, i=1,\dots ,n. \right\}.
\end{eqnarray}
\end{definition}

We equip these spaces with the following norms.
 
\begin{equation}
\|u\|_{B^{1,1/2}_{\cdot,\cdot}(Q)}= 
\|\frac{\partial^{1/2}_- u}{\partial t^{1/2}} \|_{L^2(Q)} + \|u\|_{L^p(Q)}+
\sum_{i=1}^n \|\frac{\partial u}{\partial x_i} \|_{L^p(Q)}.
\end{equation}

Computing in $\mathcal{F'}_{\cdot,\cdot}(Q)$
we see that we can represent these spaces as closed subspaces of 
the direct sum 
$L^2(Q)\oplus L^p(Q_)\oplus \cdots \oplus L^p(Q)$, and thus they are 
reflexive and separable Banach spaces in the topologies 
arising from the given norms.

Since the lateral boundary is smooth (in fact Lipschitz continuous suffices), 
we can extend an element in 
$B^{1,1/2}_{\cdot,\cdot}(Q)$ to all of ${\bf R}^n \times {\bf R}$ and then cut off in the
 space variables. 
By regularizing it is clear that functions smooth up to the boundary are dense in 
$B^{1,1/2}_{\cdot,\cdot}(Q)$. To show that
$\mathcal F_{\cdot,\cdot}(Q)$ is 
dense in $B^{1,1/2}_{\cdot,\cdot}(Q)$ we only have to prove that we can ``cut off'' in time. 
This will follow as in Lemma \ref{lem:cutoff} once we have the following result.

\begin{lemma}\label{lem:equivnorms}
If $u \in B^{1,1/2}_{\cdot,\cdot}(Q)$, then
\begin{equation}\label{eq:equivnorms}
\iiint_{\Omega \times {\bf R}\times{\bf R}} 
\left|\frac{u(x,s)-u(x,t)}{s-t}\right|^2 \, dx\,ds\,dt 
=2\pi \iint_{Q} \left|\frac{\partial_-^{1/2} u}{\partial t^{1/2}} \right|^2 \,dx\, dt.
\end{equation}
\end{lemma}

{\bf Proof.}

That $\frac{\partial_-^{1/2}u}{\partial t^{1/2}}=v$ means that
\begin{eqnarray}\nonumber
\iint_Q u(x,t)\frac{\partial_+^{1/2}\Phi(x,t)}{\partial t^{1/2}}\,dx\,dt =
\iint_Q v(x,t)\Phi(x,t)\,dx\,dt \\
;\Phi \in \mathcal{F}_{0,\cdot}(Q).
\end{eqnarray}
Now for almost every $x \in \Omega$, $\Omega \ni x \mapsto u(x,\cdot) \in L^p({\bf R})$ 
and  $\Omega \ni x \mapsto v(x,\cdot) \in L^2({\bf R})$ are well defined.
Let $S$ denote the set of common Lebesgue points. Since the Lebesgue points of a function 
can only increase by multiplication with a smooth function, by taking
limits of mean values,  
we get that
\begin{eqnarray}\nonumber
\int_{\bf R} u(x,t)\frac{\partial_+^{1/2}\Phi(x,t)}{\partial t^{1/2}}\,dt =
\int_{\bf R} v(x,t)\Phi(x,t)\,dt \\
;\Phi \in \mathcal{F}_{0,\cdot}(Q),
\end{eqnarray}
for all $x \in S$. This implies that for almost every $x\in \Omega$ the $L^p({\bf R})$ function 
$t \mapsto u(x,t)$ has half a derivative equal to $v(x,t) \in L^2({\bf R})$. 
So from the one-dimensional result it follows that
\begin{equation}
\iint_{{\bf R}\times{\bf R}} 
\left|\frac{u(x,s)-u(x,t)}{s-t}\right|^2 \,ds\,dt 
=2\pi \int_{\bf R} \left|\frac{\partial_-^{1/2} u}{\partial t^{1/2}} \right|^2 \, dt,
\end{equation}
for almost every $x\in \Omega$. Integrating with respect to $x$, the lemma follows.
$\Box$

We conclude that: 
\begin{lemma}
The space of testfunctions  
 $\mathcal{F}_{\cdot,\cdot}(Q)$ is dense in  $B^{1,1/2}_{\cdot,\cdot}(Q)$.
\end{lemma}

We now introduce the following subspace that corresponds to  zero boundary data 
on the lateral boundary $\partial \Omega \times {\bf R}$ 
and as $|t| \rightarrow \infty$.

\begin{definition}
 Let
$B^{1,1/2}_{0,\cdot}(Q)$ 
denote the closure of $\mathcal{F}_{0,\cdot}(Q)$ in the 
$B^{1,1/2}_{\cdot,\cdot}(Q)$-topology.
\end{definition}

We shall work with the following two sets of function spaces on $Q_+$.

\begin{definition}
 Let  
$B^{1,1/2}_{*,\cdot}(Q_+)$ denote the space of functions defined on $Q_+$ that can be
extended to elements in $B^{1,1/2}_{*,\cdot}(Q)$.

Furthermore  let  
$B^{1,1/2}_{*,0}(Q_+)$ denote the space of functions defined on $Q_+$ that can be
extended by zero to elements in $B^{1,1/2}_{*,\cdot}(Q)$.
\end{definition}
Here $*$ optionally stands for $\cdot$ or $0$. A zero in the first
position corresponds to zero boundary data on the lateral boundary and 
a
zero in the second position corresponds to zero initial data.

Clearly $B^{1,1/2}_{*,0}(Q_+)$ can be identified with a closed subspace of 
$B^{1,1/2}_{*,\cdot}(Q)$.

We give the following two simple lemmas concerning these spaces 
and, as in the case of the real line, we omit the easy proofs.
\begin{lemma}
  The function space $B^{1,1/2}_{*,0}(Q_+)$ becomes a 
Banach space with the norm
\begin{eqnarray}\nonumber\label{eq:bhalf0normRn}
\| u\|_{B^{1,1/2}_{*,0}(Q_+)}=\|u\|_{L^p(Q_+)} + \|\nabla_x u\|_{L^p(Q_+)}+
\left\{ \int_{Q_+} \frac{u^2(x,t)}{t} \,dt\, dx \right.\\
\left. +
\iiint_{\Omega \times {\bf R}_+ \times {\bf R}_+} 
\left(\frac{u(x,s)-u(x,t)}{s-t} \right)^2 \,dx \,ds\,dt \right\}^{1/2}.
\end{eqnarray}
\end{lemma}

\begin{lemma}
  The function space $B^{1,1/2}_{*,\cdot}(Q_+)$ becomes a 
Banach space with the norm
\begin{eqnarray}\nonumber\label{eq:bhalfnormRn}
\| u\|_{B^{1,1/2}_{\cdot,\cdot}(Q_+)}=\|u\|_{L^p(Q_+)} +\|\nabla_x u\|_{L^p(Q_+)} \\
+
\left\{
\iiint_{\Omega \times {\bf R}_+ \times {\bf R}_+} 
\left(\frac{u(x,s)-u(x,t)}{s-t} \right)^2\,dx \, ds\,dt \right\}^{1/2}.
\end{eqnarray}
Furthermore  a continuous  symmetric extension mapping 
from  $B^{1,1/2}_{*,\cdot}(Q_+)$ to  $B^{1,1/2}_{*,\cdot}(Q)$
is given by  $E_S(u)(x,t)=u(x,|t|)$.
\end{lemma}

Computing in $\mathcal F'_{\cdot,0} (Q_+)$ we can give an equivalent
characterization of  $B^{1,1/2}_{\cdot,0}(Q_+)$.

\begin{lemma}
A function $u \in L^p(Q_+)$ belongs to $B^{1,1/2}_{\cdot,0}(Q_+)$ if and only if
the $\mathcal F'_{\cdot,0} (Q_+)$-distribution derivative 
$\frac{\partial^{1/2}_+ u}{\partial t^{1/2}}$ 
belongs to $L^2(Q_+)$, and the 
$\mathcal F'_{\cdot,\cdot} (Q_+)$-distribution derivatives 
 $\nabla_x u \in L^p(Q_+)$. 
Furthermore an equivalent norm on $B^{1,1/2}_{\cdot,0}({\bf R}_+)$ 
is then given by
\begin{equation}
\|u\| = \| \nabla_x u\|_{L^p(Q_+)} +  
\| u\|_{L^p(Q_+)}+
\|\frac{\partial^{1/2}_+ u}{\partial t^{1/2}}\|_{L^2(Q_+)}.
\end{equation}
\end{lemma}

{\bf Proof.}
As on the real line.
$\Box$

Using the corresponding result on the real half-line 
and the same type of argument as in the proof of 
Lemma \ref{lem:equivnorms}, we see that an equivalent norm on 
$B^{1,1/2}_{*,\cdot}(Q_+)$  is given by 
 
\begin{equation}
\|u\| = 
\|\frac{\partial^{1/2}_- u}{\partial t^{1/2}} \|_{L^2(Q_+)} + \|u\|_{L^p(Q_+)}+
\sum_{i=1}^n \|\frac{\partial u}{\partial x_i} \|_{L^p(Q_+)},
\end{equation}
where $\frac{\partial^{1/2}_- }{\partial t^{1/2}}$ 
is understood in the $\mathcal{F'}_{\cdot,\cdot}(Q_+)$-distribution sense.

We have the following density results:
\begin{theorem}\label{th:densetestfunc}
The space of testfunctions 
 $\mathcal{F}_{\cdot,*}(Q_+)$ is dense in $B^{1,1/2}_{\cdot,*}(Q_+)$.
Furthermore the space of testfunctions 
 $\mathcal{F}_{0,*}(Q_+)$ is dense in $B^{1,1/2}_{0,*}(Q_+)$.
\end{theorem}

{\bf Proof.}
Since the boundary of $\Omega$ is smooth we have good extension operators in the 
space variables, and we can also translate the support of functions away from the lateral 
boundary without spreading the support in the time direction.
The result thus follows exactly as in Lemma \ref{lem:densereal}.
$ \Box $

We point out the following result that follows immediately from the given norms.
\begin{lemma}
The space $B^{1,1/2}_{*,0}(Q_+)$ is continuously imbedded in
$B^{1,1/2}_{*,\cdot}(Q_+)$.
\end{lemma}

We also remark that the (semi)norms 
$\|\frac{\partial^{1/2}_- u}{\partial t} \|_{L^2(Q_+)}$ 
and $\|\frac{\partial^{1/2}_+ u}{\partial t} \|_{L^2(Q_+)}$ are not equivalent. 
In fact in Lemma \ref{lem:denseimbed} below we show that 
$B^{1,1/2}_{0,0}(Q_+)$ is a dense subspace of  $B^{1,1/2}_{0,\cdot}(Q_+)$.
This is of course connected with the well known fact that if 
$u \in L^2(Q)$ and $\frac{\partial^{1/2}_- u}{\partial t^{1/2}} \in L^2(Q)$, it is 
in general impossible to define a trace on $\Omega \times \{ 0 \}$ 
(for instance the 
function $(x,t) \mapsto \log | \log |t||$ locally belongs to this space).
Still a function in $ B^{1,1/2}_{\cdot,0}(Q_+)$ is of course 
zero on $\Omega \times \{ 0\}$ 
in the sense that

\begin{equation}\label{eq:hardyineq}
 \iint_{Q_+} \frac{u^2(x,t)}{t}\,dxdt < \infty.
\end{equation}

We shall now discuss homogeneous  data on the whole parabolic boundary.

\subsection{Homogeneous data.}

We introduce the following  space of 
$\mathcal{F'}_{\cdot,\cdot}(Q)$-distributions 
 defined globally in time, but supported in $Q_+$.

\begin{definition}
Let 
\begin{equation}
B^{-1,-1/2}_{\cdot,0}(Q_+):= \left\{ \xi \in B^{1,1/2}_{0,\cdot}(Q)^*;\quad
\xi = 0 \; \mbox{in} \; \Omega \times (-\infty,0) \right\}
\end{equation}
\end{definition}

From Theorem 4.3 and Theorem 4.4 in \cite{F} follows
\begin{theorem}\label{th:homglobal}
For $T$ as defined in (\ref{eq:T}), satisfying the structural conditions (1)--(5),
\begin{equation}
T: B^{1,1/2}_{0,0}(Q_+) \longrightarrow B^{-1,-1/2}_{\cdot,0}(Q_+)
\end{equation}
is a bijection.
\end{theorem}

We shall now show that $B^{-1,-1/2}_{\cdot,0}(Q_+)$ can be identified with the
dual space of
$B^{1,1/2}_{0,\cdot}(Q_+)$.

\begin{lemma}\label{lem:bdual}
We can identify $
B^{-1,-1/2}_{\cdot,0}(Q_+)$
with
$B^{1,1/2}_{0,\cdot}(Q_+)^*$.
\end{lemma}
{\bf Remark.} Note that we here identify a subspace of 
$\mathcal{F'}_{\cdot,\cdot}(Q)$ with a subspace of 
$\mathcal{F'}_{\cdot,0}(Q_+)$.

{\bf Proof.}
Given $\xi \in B^{1,1/2}_{0,\cdot}(Q_+)^*$ we have 
(by the Hahn-Banach theorem) $u_0 \in L^2(Q_+)$ and $u_i \in L^{p'}(Q_+)$, 
$i=1,\dots,n$ such that
\begin{equation}
\langle \xi, \Phi \rangle =
\iint_{Q_+} u_0 \frac{\partial^{1/2}_- \Phi}{\partial t} +
\sum_{i=1}^n u_i \frac{\partial \Phi}{\partial x_i}\,dx dt ;\quad
\Phi \in \mathcal{F}_{0,\cdot}(Q_+).
\end{equation}
It is thus clear that we can extend this $\xi$ to all of 
$\mathcal{F}_{0,\cdot}(Q)$ by zero. Set 
\begin{equation}
\langle \xi_0, \Phi \rangle =
\iint_{Q} E_0(u_0) \frac{\partial^{1/2}_- \Phi}{\partial t} +
\sum_{i=1}^n E_0(u_i) \frac{\partial \Phi}{\partial x_i}\,dx dt ;\quad
\Phi \in \mathcal{F}_{0,\cdot}(Q),
\end{equation}
where $E_0$ denotes the operator that extends a function with $0$ to all 
of $Q$.
The mapping $B^{1,1/2}_{0,\cdot}(Q_+)^* \ni \xi \mapsto 
\xi_0 \in B^{-1,-1/2}_{\cdot,0}(Q_+)$ is clearly injective, but it is
also surjective. This follows since
given $\xi \in B^{-1,-1/2}_{\cdot,0}(Q_+)$, by  Theorem \ref{th:homglobal} 
above, there exists a (unique)
 $u_{\xi} \in B^{1,1/2}_{0,0}(Q_+)$ such that
\begin{equation}
\frac{\partial u_{\xi}}{\partial t} -\nabla_x \cdot (|\nabla_x u_{\xi}|^{p-2}
\nabla_x u_{\xi}) = \xi,
\end{equation}
i.e.
\begin{equation}\nonumber
\langle \xi, \Phi \rangle =
\iint_{Q} \frac{\partial^{1/2}_+ u_{\xi}}{\partial t}
 \frac{\partial^{1/2}_- \Phi}{\partial t}
\end{equation}
\begin{equation}
+(|\nabla_x u_{\xi}|^{p-2}
\nabla_x u_{\xi})\cdot \nabla_x \Phi  \,dx dt ;\quad
\Phi \in \mathcal{F}_{0,\cdot}(Q),
\end{equation}
and we see that $\xi$ has the required form.
$\Box$

Thus we can reformulate Theorem \ref{th:homglobal}.
\begin{theorem}\label{th:homognonlin}
For $T$ as defined in (\ref{eq:T}), satisfying the structural conditions (1)--(5),
\begin{equation}
T: B^{1,1/2}_{0,0}(Q_+) \longrightarrow B^{1,1/2}_{0,\cdot}(Q_+)^*
\end{equation}
is a bijection.
\end{theorem}
{\bf Remark.}
This theorem of course means that given $\xi \in  B^{1,1/2}_{0,\cdot}(Q_+)^*$ 
there exists a unique $u \in  B^{1,1/2}_{0,0}(Q_+)$ such that 
\begin{equation}\label{eq:tuxi}
\langle T(u), \Phi \rangle = \langle \xi , \Phi \rangle \; ; \quad 
\Phi \in  B^{1,1/2}_{0,\cdot}(Q_+).
\end{equation}
Which means precisely that
\begin{eqnarray}\nonumber
\langle \xi, \Phi \rangle =
\iint_{Q_+} \frac{\partial^{1/2}_+ u_{\xi}}{\partial t}
 \frac{\partial^{1/2}_- \Phi}{\partial t}
+ A(x,t,\nabla_x u) \cdot \nabla_x \Phi  \,dx dt \\
;\quad \Phi \in \mathcal{F}_{0,\cdot}(Q_+),
\end{eqnarray}
since  $\mathcal{F}_{0,\cdot}(Q_+)$ is dense in  $B^{1,1/2}_{0,\cdot}(Q_+)$. 

The following structure theorem for our source data space is an immediate consequence 
of the Hahn-Banach theorem.
\begin{theorem}
Given  $\xi \in  B^{1,1/2}_{0,\cdot}(Q_+)^*$ there exist functions 
$u_0 \in L^2(Q_+)$ and $u_1, \dots, u_n \in  L^{p/(p-1)}(Q_+)$ such that
\begin{equation}
\xi = \frac{\partial^{1/2}_+ u_0}{\partial t}+ \sum_{i=1}^n  \frac{\partial u_i}{\partial x_i} 
\end{equation}
in $\mathcal F'_{\cdot,0}(Q_+)$.
\end{theorem}

Our next result implies that in general it is actually enough to test our equations with  
$\mathcal{F}_{0,0}(Q_+)$
instead of  $\mathcal{F}_{0,\cdot}(Q_+)$.

\begin{lemma}\label{lem:denseimbed}
The continuous imbedding
\begin{equation}
B^{1,1/2}_{0,0}(Q_+) \hookrightarrow B^{1,1/2}_{0,\cdot}(Q_+)
\end{equation}
is dense.
\end{lemma}
{\bf Proof.}
It is enough to show that 
 if $\xi \in  B^{1,1/2}_{0,\cdot}(Q_+)^*$ and
$\langle \xi,\Phi \rangle =0$ for all $ \Phi \in B^{1,1/2}_{0,0}(Q_+)$, then 
$\xi =0$.

Now given  $\xi \in  B^{1,1/2}_{0,\cdot}(Q_+)^*$, by Theorem \ref{th:homognonlin}, 
there exists a unique $u_{\xi} \in  B^{1,1/2}_{0,0}(Q_+)$ such that
\begin{equation}
\frac{\partial u_{\xi}}{\partial t} - \nabla_x \cdot(|\nabla_x u_{\xi}|^{p-2}
\nabla_x u_{\xi}) = \xi.
\end{equation}
Now if $\langle \xi, \Phi \rangle =0$ for all $\Phi \in B^{1,1/2}_{0,0}(Q_+)$, 
then with $\Phi=u_{\xi}$ we get
\begin{equation}
\iint_{Q_+}|\nabla_x u_{\xi}|^p \,dx\,dt =0.
\end{equation}
By the Poincar\'e inequality $u_{\xi}=0$, and so $\xi=0$.
$\Box$

\subsection{Non-homogeneous initial data.}

We will first introduce the space that will carry the initial data. 
In the definition, all derivatives should be 
understood in the $\mathcal{F'}_{\cdot,\cdot}(Q_+)$-distribution sense.
\begin{definition}
Let 
\begin{eqnarray}\nonumber
 B_I(Q_+)= \left\{ u \in B^{1,1/2}_{0,\cdot}(Q_+)\cap 
C_b([0,\infty),L^2(\Omega)) \right.\\
\left. ;\frac{\partial u}{\partial t} 
\in L^{p'}({\bf R}_+, W^{-1,p'}(\Omega)) \right\} .
\end{eqnarray}
\end{definition}
Here $C_b([0,\infty),L^2(\Omega))$ denotes the space of 
bounded continuous functions 
from $[0,\infty)$ into $L^2(\Omega)$, and  
$\frac{\partial u}{\partial t} 
\in L^{p'}({\bf R}_+, W^{-1,p'}(\Omega))$ means exactly that
\begin{equation}
|\langle u, \frac{\partial \Phi}{\partial t} \rangle | 
\leq C \|\nabla_x \Phi \|_{L^p(Q_+)}\quad; \quad \Phi \in \mathcal{F}_{0,0}(Q_+),
\end{equation}
for some constant $C>0$. The smallest possible constant is by definition 
$\|\frac{\partial u}{\partial t}\|_{L^{p'}({\bf R}_+, W^{-1,p'}(\Omega))}$.

We equip $B_I(Q_+)$ with the following norm 
\begin{eqnarray}\nonumber
\|u\|_{B_I(Q_+)}:= \|u\|_{B^{1,1/2}_{0,\cdot}(Q_+)}+
\sup_{t \in {\bf R}_+} \| u(\cdot,t) \|_{L^2(\Omega)} \\
+ \|\frac{\partial u}{\partial t}\|_{L^{p'}({\bf R}_+, W^{-1,p'}(\Omega))}.
\end{eqnarray}

Using Theorem \ref{th:homognonlin} and the monotonicity of $A(x,t,\cdot)$ we shall 
now prove that we always have a unique solution in $B_I(Q_+)$ to the following  
initial value problem.

\begin{theorem}
Given $u_0 \in L^2(\Omega)$, there exists a unique element  
$u \in B_I(Q_+)$ such that
\begin{subequations}\label{eq:initialvalues}
\begin{align}
 \frac{\partial u}{\partial t}  - \nabla_x\cdot  A(x,t,\nabla_x u)
 &=0 \quad
\mbox{in $\mathcal{F'}_{\cdot,\cdot}(Q_+)$}\\
u&=u_0 \quad \mbox{on $ \Omega \times \{0\} $}.
\end{align}
\end{subequations}
\end{theorem}

{\bf Proof.}

Uniqueness follows immediately from the monotonicity of $A(x,t,\cdot)$ by pairing with
a cut off function in time multiplied with the difference of two solutions.
To prove existence 
we first note that if $u_0 \in \mathcal D (\Omega)$, we can extend it for instance 
to a smooth testfunction $U_0 \in \mathcal D (\Omega \times (-2,2))$ 
such that $U_0(x,t)=u_0(x)$ when $-1<t<1$.

Since $\frac{\partial U_0}{\partial t} \in  B^{1,1/2}_{0,\cdot}(Q_+)^*$, 
by Theorem \ref{th:homognonlin}, we know that there exists a unique 
$w \in B^{1,1/2}_{0,0}(Q_+)$ such that
\begin{equation}
 \frac{\partial w}{\partial t}  - \nabla_x\cdot  A(x,t,\nabla_x w + \nabla_x U_0)
 =-\frac{\partial U_0}{\partial t} \quad
\mbox{in $ B^{1,1/2}_{0,\cdot}(Q_+)^*$}.
\end{equation}
Then clearly $u= (w+U_0) \in  B^{1,1/2}_{0,\cdot}(Q_+)$ solves
(\ref{eq:initialvalues}),
and the initial value 
is taken in the 
sense that
\begin{equation}
\iint_{\Omega \times (0,1)} \frac{(u(x,t)-u_0(x))^2}{t}\,dx\,dt < \infty.
\end{equation}
By standard arguments it follows from (\ref{eq:initialvalues}) 
that $u \in B_I(Q_+)$ and so the 
initial data is actually taken in $C_b([0,\infty),L^2(\Omega))$-sense.

Given $u_0 \in L^2(\Omega)$ we now choose a sequence 
$\mathcal D(\Omega) \ni u^n_0 \longrightarrow u_0$ in $L^2(\Omega)$.

Let $u^n$ denote the solution of (\ref{eq:initialvalues}) with initial data $u^n_0$.
By testing with $u^n \chi$, where $\chi$ is a standard cut off function in time, 
 in (\ref{eq:initialvalues}), we get that 
\begin{equation}
\sup_{t \in {\bf R}_+} \int_{\Omega} (u^n -u^m)^2(x,t)\,dx \leq 
\int_{\Omega} (u^n_0 -u^m_0)^2(x)\,dx.
\end{equation}
It is also clear that $\|\nabla_x u^n \|_{L^p(Q_+)}$ is bounded by a 
constant independent of $n$.

Finally we note that we can extend $u^n$ symmetrically to $Q$ and the 
extended function $E_S(u^n) \in  B^{1,1/2}_{0,\cdot}(Q)$ will  satisfy
$\frac{\partial E_S(u^n)}{\partial t} 
\in L^{p'}({\bf R}, W^{-1,p'}(\Omega))$. 

We then have
\begin{equation}
\iint_Q \frac{\partial_-^{1/2} E_S(u^n)}{\partial t^{1/2}}  
\frac{\partial_-^{1/2} \Phi_k}{\partial t^{1/2}} \,dx\,dt
= \int_{\bf R}  \langle 
\frac{\partial E_S(u^n)}{\partial t}, h(\Phi_k) \rangle \,dt, 
\end{equation}
for a sequence $\mathcal {F}_{0,\cdot}(Q) \ni \Phi_k \rightarrow E_S(u^n)$ 
in  $B^{1,1/2}_{0,\cdot}(Q)$.

This implies that 
$\|\frac{\partial_-^{1/2} E_S(u^n)}{\partial t^{1/2}}\|_{L^2(Q_+)}$ 
is bounded by a constant independent of $n$. 

We conclude that
$\| u^n \|_{B_I(Q_+)} \leq C$, where $C <\infty $ is a constant independent of $n$.
 
We can now extract a weakly convergent subsequence and in fact, as we have seen,
 we actually have strong convergence in  $C_b([0,\infty),L^2(\Omega))$ and thus 
the limit function satisfies the initial conditions.

Finally a  Minty argument using the monotonicity of  $A(x,t,\cdot)$ shows that the 
limit function solves (\ref{eq:initialvalues}). The theorem follows.
$\Box$

\subsection{Fully non-homogeneous initial-boundary values.}

We shall now introduce the function space that will 
carry both initial and lateral boundary data.

Since we have continuous imbeddings
$B^{1,1/2}_{0,\cdot}(Q_+)\hookrightarrow B^{1,1/2}_{\cdot,\cdot}(Q_+)$ and \\
$B_I(Q_+) \hookrightarrow B^{1,1/2}_{\cdot,\cdot}(Q_+)$, the following
definition makes sense.
\begin{definition}\label{def:xspace}
Let
\begin{equation}
X^{1,1/2}(Q_+)=  B^{1,1/2}_{\cdot,0}(Q_+) + B_I(Q_+),
\end{equation}
be equipped with the norm
\begin{equation}
\|u \|_{X^{1,1/2}(Q_+)} = \inf_{(u_1,u_2) \in K_u} \left(
\|u_1\|_{B^{1,1/2}_{\cdot,0}(Q_+)} + \| u_2\|_{B_I(Q_+)}\right),
\end{equation}
where the infimum is taken over the set
\begin{equation}
K_u=\left\{ (u_1,u_2);\; u_1+u_2 =u,\; u_1 \in  B^{1,1/2}_{\cdot,0}(Q_+),\;
u_2 \in B_I(Q_+) \right\}.
\end{equation}
\end{definition}

The following imbeddings are immediate
\begin{eqnarray}
\| u\|_{X^{1,1/2}(Q_+)} \leq \|u\|_{B_{\cdot,0}^{1,1/2}(Q_+)};\quad 
 u \in B_{\cdot,0}^{1,1/2}(Q_+),\\
\| u\|_{X^{1,1/2}(Q_+)} \leq \|u\|_{B_I(Q_+)};\quad 
 u \in B_I(Q_+),\\
\| u\|_{B^{1,1/2}_{\cdot,\cdot}(Q_+)} \leq C \| u\|_{X^{1,1/2}(Q_+)};\quad
u \in X^{1,1/2}(Q_+).
\end{eqnarray}

For an element in $X^{1,1/2}(Q_+)$ 
we can always define the trace on $\Omega \times \{0\}$.
\begin{theorem}\label{th:tracex12}
There exists a continuous linear and surjective trace operator
\begin{equation}
Tr_0:X^{1,1/2}(Q_+) \longrightarrow L^2(\Omega).
\end{equation}
There also exists a bounded extension operator
\begin{equation}
E : L^2(\Omega) \longrightarrow X^{1,1/2}(Q_+)
\end{equation}
such that $Tr_0 \circ E = Id_{L^2(\Omega)}$.
\end{theorem}
{\bf Proof.} Given $u \in X^{1,1/2}(Q_+)$, there exist 
$u_1 \in B_{\cdot,0}^{1,1/2}(Q_+)$ and $u_2 \in B_I(Q_+)$ such that
$u=u_1+u_2$. Since $ u_2 \in B_I(Q_+) \Longrightarrow u_2 \in 
C_b([0,+\infty),L^2(\Omega))$, $u_2|_{\Omega \times \{0\}}$ is a well defined
element of $L^2(\Omega)$. We now define $u|_{\Omega \times \{0\}}=
u_2|_{\Omega \times \{0\}}$. We have to show that this is independent of 
the decomposition of $u$, but if we have two different decompositions
$u_1+u_2 = v_1+v_2$ as above, then 
$(u_2-v_2) \in B_I(Q_+) \cap B_{\cdot,0}^{1,1/2}(Q_+)$, which implies that
\begin{equation}
\iint_{\Omega \times (0,+\infty)} \frac{(u_2-v_2)^2(x,t)}{t}\,dxdt < +\infty,
\end{equation}
and so $u_2(\cdot,0)=v_2(\cdot,0)$ since they both belong to 
$C_b([0,+\infty),L^2(\Omega))$.

Now  
\begin{equation}
\|u(\cdot,0)\|_{L^2(\Omega)} = \| u_2(\cdot,0)\|_{L^2(\Omega)}
\leq C \|u_2 \|_{B_I(Q_+)},
\end{equation}
for any decomposition $u=u_1+u_2$ as above. Taking the infimum over all such
decompositions gives:
\begin{equation}
\| u(\cdot,0)\|_{L^2(\Omega)} \leq C \|u\|_{X^{1,1/2}(Q_+)},\; 
u \in X^{1,1/2}(Q_+).
\end{equation}

Now
given $u_0 \in L^2(\Omega)$, let $E(u_0)$ be the (unique) solution in
$B_I(Q_+)$ of the initial value problem:
\begin{subequations}
\begin{align}
 \frac{\partial u}{\partial t} - \nabla_x \cdot (|\nabla_x u|^{p-2}\nabla_x u)
 &=0 \quad
\mbox{in $Q_+=\Omega \times {\bf R_+}$}\\
u&=u_0 \quad \mbox{on $ \Omega \times \{0\} $}.
\end{align}
\end{subequations}
Clearly this extension map satisfies 
$Tr_0 \circ E = Id_{L^2(\Omega)}$ and furthermore
\begin{equation}
\|E (u_0)\|_{B_I(Q_+)} \leq C \|u_0\|_{L^2(\Omega)},
\end{equation}
and thus
\begin{equation}
\|E (u_0)\|_{X^{1,1/2}(Q_+)} \leq C \|u_0\|_{L^2(\Omega)}.
\end{equation}
$\Box$

{\bf Remark.} Note that if $p=2$ the extension map is linear.

\begin{theorem}
We have the following imbedding:
\begin{equation}
\| u\|_{B^{1,1/2}_{\cdot,0}(Q_+)} \leq C \| u\|_{X^{1,1/2}(Q_+)}; \; u
\in B^{1,1/2}_{\cdot,0}(Q_+).
\end{equation}
\end{theorem}
{\bf Proof.} If $u \in B^{1,1/2}_{\cdot,0}(Q_+)$, and 
$u=u_1+u_2$ with $u_1 \in B_{\cdot,0}^{1,1/2}(Q_+)$ and $u_2 \in B_I(Q_+)$, 
then $u_2(\cdot,0)=0$ since $u_2 \in  B_{\cdot,0}^{1,1/2}(Q_+) \cap 
 B_I(Q_+)$. Thus $u_2$ can be extended by zero to all of $Q$.
 Since, by a continuity argument, 
\begin{equation}
\| \frac{\partial_+^{1/2} u_2}{\partial t} \|_{L^2(Q_+)}^2
= -\int_{{\bf R}_+} \langle \frac{\partial u_2}{\partial t}, h(u_2)
\rangle \,dt, \quad
u_2 \in  B^{1,1/2}_{\cdot,0}(Q_+) \cap B_I(Q_+).
\end{equation}
We get
$$
\| u_1\|_{ B_{\cdot,0}^{1,1/2}(Q_+)} + \|u_2\|_{ B_I(Q_+)}
$$
$$
\geq C\left( \| u_1\|_{ B_{\cdot,0}^{1,1/2}(Q_+)} + 
\|u_2\|_{  B_{\cdot,0}^{1,1/2}(Q_+)} \right)
$$
\begin{equation}
\geq C\|u_1+u_2\|_{ B_{\cdot,0}^{1,1/2}(Q_+)}= 
C\|u \|_{ B_{\cdot,0}^{1,1/2}(Q_+)},
\end{equation}
where $C>0$.
Taking the infimum concludes the proof.
$\Box$
We immediately get the following
\begin{corollary}
There exist constants $C_1,C_2>0$ such that
\begin{equation}
C_1 \| u \|_{ B_{0,0}^{1,1/2}(Q_+)} \leq  \| u \|_{ X^{1,1/2}(Q_+)} 
\leq  C_2 \| u \|_{ B_{0,0}^{1,1/2}(Q_+)};\; u \in  
 \mathcal F_{0,0}(Q_+).
\end{equation}
Thus  $B_{0,0}^{1,1/2}(Q_+)$ is the closure of 
$\mathcal F_{0,0}(Q_+)$ in the $X^{1,1/2}(Q_+)$-norm topology.
\end{corollary}

We are now ready to state our main theorem.
\begin{theorem}\label{th:mainnonlinear}
Given $f \in  B_{0,\cdot}^{1,1/2}(Q_+)^*$ and 
$g \in X^{1,1/2}(Q_+)$, there exists a unique element 
$u \in  X^{1,1/2}(Q_+)$ such that
\begin{subequations}\label{eq:proofnonlinear}
\begin{align}
 \frac{\partial u}{\partial t} - \nabla_x \cdot (A(x,t,\nabla_x u))
 &=f \quad
\mbox{in}\; \mathcal F'_{\cdot,\cdot}(Q_+) \\
u-g \in & \,  B_{0,0}^{1,1/2}(Q_+).
\end{align}
\end{subequations}
\end{theorem}
{\bf Proof.}

Let $w=u-g$. Then (\ref{eq:proofnonlinear}) is equivalent to
\begin{subequations}\label{eq:proofhomnonlinear}
\begin{align}
 \frac{\partial w}{\partial t} - \nabla_x \cdot (A(x,t,\nabla_x (w+g)))
 &=f-\frac{\partial g}{\partial t} \; \mbox{in}\;   \mathcal F'_{\cdot,\cdot}(Q_+)\\
w \in & \, B_{0,0}^{1,1/2}(Q_+).
\end{align}
\end{subequations}
Here  
$\frac{\partial g}{\partial t} \in  \mathcal{F'}_{\cdot,\cdot}(Q_+)$ 
 has a unique  extension to an element in $B_{0,\cdot}^{1,1/2}(Q_+)^*$.
 In fact, if 
$g \in  X^{1,1/2}(Q_+)$, we can write $g=g_1+g_2$, where 
$g_1 \in  B_{\cdot,0}^{1,1/2}(Q_+)$ and $g_2 \in  B_I(Q_+)$. Thus
\begin{equation}\nonumber
|\langle g, \frac{\partial \Phi}{\partial t} \rangle | =
|\langle g_1, \frac{\partial \Phi}{\partial t} \rangle +
\langle g_2, \frac{\partial \Phi}{\partial t} \rangle |
\end{equation}
\begin{equation}
\leq C \left( \|g_1\|_{ B_{\cdot,0}^{1,1/2}(Q_+)} +
\| g_2\|_{B_I(Q_+)} \right) \| \Phi \|_{B_{0,\cdot}^{1,1/2}(Q_+)};\; 
\Phi \in \mathcal{F}_{0,0}(Q_+).
\end{equation}
 Since, by Lemma \ref{lem:denseimbed} and Theorem \ref{th:densetestfunc}, 
$\mathcal F_{0,0}(Q_+)$ is dense in $B_{0,\cdot}^{1,1/2}(Q_+)$, 
it is clear that we have a unique extension.
If the function $A(\cdot,\cdot,\cdot)$ satisfies the structural 
conditions 1--5 given above, then
also $A(\cdot,\cdot,\cdot + g)$, with $g \in  X^{1,1/2}(Q_+)$, 
 satisfies the same structural conditions 
(with new constants $\lambda , \Lambda$ and functions $H,h$ depending on 
$g$). Thus Theorem \ref{th:homognonlin}, and the remark following 
Theorem \ref{th:homognonlin}, tell us that 
(\ref{eq:proofhomnonlinear}) has a unique solution.
This implies that $u=w+g$ is the unique solution to (\ref{eq:proofnonlinear}).
$\Box$

{\bf Remark.}  Note that since $\mathcal D(Q_+)$ is densely continuously imbedded 
in  $\mathcal F_{0,0}(Q_+)$ it is equivalent to demand that (\ref{eq:proofnonlinear}) 
should hold in  $\mathcal D'(Q_+)$. 

We shall conclude with a comment on the linear case. 

The function spaces we have 
introduced so far coincides with well known function spaces 
existing in the literature when $p=2$. When $p=2$ we shall follow 
existing notation and replace $B$ with $H$ for all spaces (for instance if $p=2$ 
we shall write $H^{1,1/2}_{0,\cdot}(Q_+)$ instead of
$B^{1,1/2}_{0,\cdot}(Q_+)$ and so on).

The Sobolev space 
$H_{\cdot,\cdot}^{1/2,1/4}(\partial \Omega \times {\bf R}_+)$ below is 
 defined by pull-backs in local charts on $\partial \Omega$.

\begin{theorem}\label{th:lateraltrace}
If $p=2$ there exists a linear,
 continuous and surjective trace operator
\begin{equation}
Tr: X^{1,1/2}(Q_+) \longrightarrow 
H_{\cdot,\cdot}^{1/2,1/4}(\partial \Omega \times
{\bf R}_+).
\end{equation}
There also exists a continuous and linear extension operator
\begin{equation}
E: H_{\cdot,\cdot}^{1/2,1/4}(\partial \Omega \times {\bf R}_+) \longrightarrow 
 X^{1,1/2}(Q_+),
\end{equation}
such that $Tr \circ E =
Id|_{H_{\cdot,\cdot}^{1/2,1/4}(\partial \Omega \times {\bf R})} $.
\end{theorem}
{\bf Proof.}
Using a partition of unity argument and  the Fourier multiplier operators
\begin{equation}
m_s(D)u = ((1+i2\pi\tau +4\pi^2|\xi|^2)^{-s}\hat{u})^{\vee};\quad s\in {\bf R},
\end{equation}
which preserves forward support in time, and have the property that
\begin{equation}
m_s(D)\left(L^2({\bf R}^n \times {\bf R})\right)= 
H_{\cdot,\cdot}^{2s,s}({\bf R}^n \times {\bf R}),
\end{equation}  
we can construct continuous linear operators:
\begin{equation}
Tr: H^{1,1/2}_{\cdot,0}(Q_+) \longrightarrow 
H_{\cdot,\cdot}^{1/2,1/4}(\partial \Omega \times
{\bf R}_+)
\end{equation}
and 
\begin{equation}
E: H_{\cdot,\cdot}^{1/2,1/4}(\partial \Omega \times {\bf R}_+) \longrightarrow 
 H^{1,1/2}_{\cdot,0}(Q_+),
\end{equation}
such that  $Tr \circ E =Id|_{H_{\cdot,\cdot}^{1/2,1/4}
(\partial \Omega \times {\bf R})} $.
Now given $u \in X^{1,1/2}(Q_+)$, let $u=u_1+u_2$ where $u_1 \in 
 H^{1,1/2}_{\cdot,0}(Q_+)$ and $u_2 \in H_I(Q_+)$. We define 
$u|_{\partial \Omega \times {\bf R}_+} = 
u_1|_{\partial \Omega \times {\bf R}_+}$. This definition is 
independent of the decomposition of $u$. In fact, if $u_1+u_2 =v_1+v_2$ are
two decompositions as above, then $u_1-v_1 \in L^2({\bf R}_+, H_0^1(\Omega))$, 
and so $(u_1-v_1)|_{\partial \Omega \times {\bf R}}=0$.

Now 
\begin{equation}
\|Tr(u)\|_{H_{\cdot,\cdot}^{1/2,1/4}(\partial \Omega \times {\bf R}_+)} \leq 
C\| u_1\|_{ H^{1,1/2}_{\cdot,0}(Q_+)},
\end{equation}
for any decomposition. Taking the infimum proves the continuity of $Tr$.
The continuity of the extension operator $E$ follows from the imbedding
 $H_{\cdot,0}^{1,1/2}(Q_+) \hookrightarrow X^{1,1/2}(Q_+)$.
$\Box$

Combining our trace theorems with Theorem \ref{th:mainnonlinear} gives us in
the linear case:
\begin{theorem}\label{th:mainlinear}
If
\begin{equation}
Tu = \frac{\partial u}{\partial t} -\nabla_x \cdot (A(x,t,\nabla_x u)),
\end{equation}
is a linear operator, satisfying the structural conditions 1--5 above,
then
\begin{equation}\nonumber
 X^{1,1/2}(Q_+) \ni u \mapsto (Tu, u|_{\partial \Omega \times {\bf R}_+},
u|_{\Omega \times \{0\}})
\end{equation}
\begin{equation}
\in  H^{1,1/2}_{0,\cdot}(Q_+)^* \times  
H_{\cdot,\cdot}^{1/2,1/4}(\partial \Omega \times {\bf R}_+)\times L^2(\Omega),
\end{equation}
is a linear isomorphism.
\end{theorem}
\vspace{2cm}
\begin{center}
{\bf ACKNOWLEDGEMENTS.}
\end{center}

  I thank Johan R{\aa}de for useful remarks and stimulating discussions
  in connection with this work,
  and Anders Holst, Per-Anders Ivert and Stefan Jakobsson for reading
  and commenting on this paper.

\end{document}